\documentclass[12pt]{amsart}
\usepackage{amsmath}
\usepackage{geometry,amsfonts,amssymb,amsthm,txfonts,pxfonts,amscd} 
\def\struckint{\mathop{%
\def\mathpalette##1##2{\mathchoice{##1\displaystyle##2}%
  {##1\textstyle##2}{##1\scriptstyle##2}{##1\scriptscriptstyle##2}}%
\mathpalette
{\vbox\bgroup\baselineskip0pt\lineskiplimit-1000pt\lineskip-1000pt
\halign\bgroup\hfill$}
{##$\hfill\cr{\intop}\cr\diagup\cr\egroup\egroup}%
}\limits}
\newtheorem{theorem}{Theorem}

\theoremstyle{remark}

\DeclareMathOperator{\tr}{tr}

\begin{document}

\title{Golden-Thompson from Davis}

\author{Igor Rivin}
%\address{Department of Mathematics, Temple University, Philadelphia}
\address{School of Mathematics, Institute for Advanced Study, Princeton}
\email{rivin@ias.edu}
\thanks{Igor Rivin was supported by the National Science Foundation. He would also like to thank the Institute for Advanced Study for its hospitality and Nikhil Srivastava for interesting conversations.}
\date{\today} % delete this line to display the current date

%%% BEGIN DOCUMENT

\begin{abstract}
We give a very short proof of the Golden-Thompson inequality
\end{abstract}
\maketitle
%\tableofcontents

The Golden-Thompson inequality (due independently to Golden \cite{goldenineq} and Thompson \cite{thompsonineq}) states that for two hermitian matrices $A$ and $B,$ the following inequality holds:
\begin{equation}
\label{gtineq}
\tr \exp(A+B) \leq \tr \exp A \tr \exp B.
\end{equation}
The Golden-Thompson inequality (which is widely used in statistical physics and the theory of random matrices) is usually proved in a rather painful way using the power series of the matrix exponential function. In this note we point out that the Golden-Thompson follows immediately from the celebrated (in some circles) and underappreciated (in other circles) theorem of Chandler Davis, which states:
\begin{theorem}[Chandler Davis' Theorem] Let $f$ be a unitarily invariant function on the set of hermitian matrices. Then $f$ is convex if and only if its restriction to the set of \emph{diagonal} matrices is a convex symmetric function.
\end{theorem}
A function is \emph{unitarily invariant} if $f(U^t A U) = f(A)$ for any unitary $U.$
The original proof (\cite{davispap}) is very short and elegant, and the reader is advised to read Davis' (one page) paper. A different two-page proof has been given by the author (\cite{mydavis}  Both proofs go through verbatim if we replace "hermitian" by "symmetric" and "unitary" by "orthogonal".

Now, the Golden-Thompson inequality states that the function $F(x)=\log \tr \exp (x)$ is convex on the set of hermitian matrices, and since it is obviously unitarily invariant, we need only check convexity of the function 
$f(\lambda_1,  \dots, \lambda_n) = \log(\sum_{i=1}^n \exp(\lambda_i)).$ Now, for the function $f$ we can compute the hessian without too much trouble: 
\begin{gather}
H_{ii}=\dfrac{\partial^2 f}{\partial \lambda_i^2} = \dfrac{\sum_{j\neq i} \exp(x_i+x_j)}{(\sum_{k=1}^n \exp(x_k))^2},\\
H_{ij} = \dfrac{\partial^2 f}{\partial \lambda_i \partial \lambda_j} = 
-\dfrac{\exp(x_i+x_j)}{(\sum_{k=1}^n\exp(x_k))^2}
\end{gather}
It is easy to see that $H= D K D,$ where D is the diagonal matrix with the $i$-th entry equal to 
$\exp(x_i)/(\sum_j \exp(x_j),$ and $K$ is the Laplacian matrix of the complete graph on $n$ vertices (in other words, $K_{ii}=n-1,$ and $K_{ij}=-1$ when $i\neq j.$) Since $K$ is a connected graph, its laplacian matrix is positive semidefinite (with dimension of the null-space equal to 1), the Golden-Thompson inequality follows immediately.

\bibliographystyle{plain}
\bibliography{gs}
\end{document}